\documentclass[11pt]{amsart}
\usepackage{graphicx}
\usepackage{amssymb}
\usepackage{epstopdf}
\usepackage{hyperref}

\def\R{\mathbf{R}}

\newtheorem{innercustomgeneric}{\customgenericname}
\providecommand{\customgenericname}{}

\theoremstyle{plain}
  \newtheorem{thm}{Theorem}[section]

\begin{document}

\title[Central Configurations]{``Twist Vectors" for Central Configuration Equations}
\author{Marshall Hampton}

\begin{abstract}
A new coordinate system on the tangent space to planar configurations is introduced to simplify some calculations on central configurations and relative equilibria in the $N$-body problem with a homogeneous potential, which includes the case of Newtonian gravity.  These coordinates are applied to some problems on four-body central configurations to illustrate their utility.
\end{abstract}

\maketitle

\section{Introduction}

In this paper we consider the generalization of relative equilibria to `central configurations' for a homogeneous family of potentials that includes the Newtonian case.   As the classical Newtonian $N$-body problems do not possess static equilibria, the relative equilibria (equilibria in a uniformly rotating coordinate system) are important anchors from which we can extend our understanding of the dynamics of these systems.  We introduce some coordinates on the tangent space of configurations that we believe provides a simplifying context for some known equations and results, and which will hopefully allow further progress.

There is a large literature on central configurations, from which we recommend the book chapter by Moeckel \cite{LMS} for a wide ranging summary and background.  

Most of the natural questions concerning central configurations in the three-body problem have been answered, such as their existence, geometry, and number \cite{Euler_col, Lagrange_3}, and their linear stability \cite{gascheau1843, routh1875, roberts02}.  However for the four-body and higher $N$-body problems many open questions remain \cite{albouy2012some}, including the most famous problem of the finiteness of planar central configurations \cite{chazy18,wintner41,smale98}.

\section{Twist Vectors and Central Configurations}

We define a (scale-normalized) central configuration by the conditions
\begin{equation} \label{defcc}
(q_i - c) - \sum_{j = 1, j \neq i}^n \frac{m_j (q_i - q_j)}{r_{i,j}^A} = 0
\end{equation}
where $q_i \in \R^d$ is the position of particle $i$, and $r_{i,j} = |q_i - q_j|$ are the mutual distances between the particles.  $A>2$ is a real parameter.  The vector $c$ is the center of mass, 
$$c = \frac{1}{M} \sum_{i=1}^n m_i q_i,$$
with $M = \sum_{i=1}^{n} m_i$ the total mass.  

We will focus on the most studied case of planar configurations $(d=2)$. When considering the space of configurations we will use the flattened coordinates in $\R^{2n}$: $(x_1, y_1, x_2, y_2, \ldots, x_n, y_n)$.

An equivalent, sometimes more convenient version of the equations \ref{defcc} is

\begin{equation} \label{defcc2}
\sum_{k = 1, j \neq i}^n m_k (q_i - q_k)(r_{i,k}^{-A} - 1) = 0 .
\end{equation}
These equations can be viewed as critical point equations in several different ways; we choose to consider them as requiring the gradient of $f$ to be zero where 

$$ f = MI/2 + U/(A-2) $$
in which $I = \sum_{i=1}^{n} m_i |q_i - c|^2$ is the moment of inertia around the center of mass, and $U = \sum_{i=1}^{n} m_i m_j r_{i,j}^{-A+2}$ is a homogeneous potential function.  For the special case $A=2$, which is of some interest in the theory of fluid vorticity \cite{HelmholtzV,kirchhoff_vorlesungen_1883,aref1992grobli}, we can instead use $U = \sum_{i=1}^{n} m_i m_j \ln(r_{i,j})$ (without the division by $(A-2)$).

At a critical point Euler's relation for homogeneous functions implies that 

$$U = M I. $$

In terms of mutual distances this condition can be written as

$$ \sum_{i>j}  m_i m_j r_{i,j}^2 \left ( r_{i,j}^{-A} -  1 \right )  = 0 . $$

Our choice of definitions means we have set the scale for the overall size of the central configurations; because of the rotational and translation invariance of the Newtonian dynamics the main point of interest is the shape of the central configurations (i.e. modulo uniform rotation and translation).  Therefore it is helpful to recast the equations (\ref{defcc2}) in terms of translation- and rotation-invariant geometric quantities such as the mutual distances $r_{i,j}$.  One way to do that is to take the inner product of each (vector-valued) equation in  (\ref{defcc2}) with the difference vector $q_i - q_j$ for some choice of $j \neq i$:

\begin{equation} \label{elem_aceqs}
b_{i,j} = \sum_{k=1, k \neq i}^n m_k (r_{i,k}^{-A} - 1) (r_{i,j}^2 + r_{i,k}^2 - r_{j,k}^2) = 0 ,
\end{equation}
which we call the `asymmetric Albouy-Chenciner equations' (originally coined by Gareth Roberts).  When summed in symmetric pairs, $g_{i,j} = b_{i,j} + b_{j,i}$, we obtain the (symmetric) Albouy-Chenciner equations  \cite{albouy_probleme_1997, hampton_finiteness_2005}.  

The fact that the Albouy-Chenciner equations involve only the mutual distances can be both an advantage or disadvantage.  For the planar case, once $n>3$ the mutual distances are no longer independent, and it can be a major complicating factor to add in their constraints (e.g. using Cayley-Menger determinants).  

We will focus on an alternative, which has been described many times before, sometimes called the Laura-Andoyer equations \cite{laura1905sulle, andoyer1906main, meyer1933solutions, hagihara1970celestial}, and extend them to a framework that also includes the Hessian of $f$ in order to exploit Morse theory \cite{milnor2016morse, bott1988morse}.    Morse theory can be used to obtain lower bounds on the number of central configurations, and the Hessian index is also related to the stability of the central configurations \cite{SunHu09, Baru14,Montaldi2015,hampton_homog2019}.

In order to make our formulas more concise, we will use the notation $x_{i,k} = x_i - x_k$, $y_{i,j} = y_i - y_j$, $q_{i,j} = (x_{i,j}, y_{i,j})$, $R_{i,k} = r_{i,k}^{-A}$, $S_{i,k} = r_{i,k}^{-A} - 1$, $T_{i,k} = r_{i,k}^{-A-2}$,  $A_{i,j,k} = r_{i,j}^2 + r_{i,k}^2 - r_{j,k}^2$, and 

$$ \Lambda_{i,j,k} = q_{i,j} \wedge q_{i,k} = ( x_{i,j} y_{i,k} - x_{i,k} y_{i,j}) ,$$
which is twice the signed area of the points $q_i$, $q_j$, and $q_k$.

The $\Lambda_{i,j,k}$ change sign under odd permutations of the indices, also note $A_{i,j,k} = A_{i,k,j}$ and $A_{i,j,j} = 2 r_{i,j}^2$. 

With this notation the equations (\ref{elem_aceqs}) become simply

$$ b_{i,j} = \sum_{k=1, k \neq i}^n m_k S_{i,k} A_{i,j,k} = 0 . $$

Note that in the sum for $b_{i,j}$, the term with $k=j$ plays a somewhat special role, so it can be convenient to separate it out:

$$ b_{i,j} = 2 m_j S_{i,j} + \sum_{k \notin \{i,j\}} m_k S_{i,k} A_{i,j,k} = 0 .$$

Our main new tool is to introduce what we will call `twist vectors'.  These configuration space tangent vectors rotate a pair of points around their center of mass.  The nonzero components of the twist vectors at a configuration of points with coordinate $ q_i = (x_i, y_i)$ are 

$$ \left ( v_{i,j} \right )_{x_i} = m_j y_{i,j}, \ \ \ \ \left ( v_{i,j} \right )_{y_i} = - m_j x_{i,j} , $$ 

$$ \left ( v_{i,j} \right )_{x_j} = - m_i y_{i,j}, \ \ \ \ \left ( v_{i,j} \right )_{y_j} = m_i x_{i,j} . $$ 

The vectors $v_{i,j}$ are all perpendicular to the vectors $\nabla c_x = (m_1, 0, m_2, 0, \ldots, m_n, 0)/M$ and $\nabla c_y = (0, m_1, 0, m_2, \ldots, 0, m_n)/M$, so the orbit of a configuration under the flow generated by any linear combination of the $v_{i,j}$ will have a constant center of mass $c = (c_x, c_y)$.

The $v_{i,j}$ are also all tangent to the level surfaces of the moment of inertia, $I$, since $\nabla I = 2(m_1 x_1, m_1 y_1, \ldots, m_n x_n, m_n y_n)$ and thus

$$ \nabla I \cdot v_{i,j} = 2 m_i m_j (y_{i,j} x_i - x_{i,j} y_i - y_{i,j} x_j + x_{i,j} y_j) = 0 . $$

\begin{center}
\begin{figure} 
\includegraphics[width=2.5in]{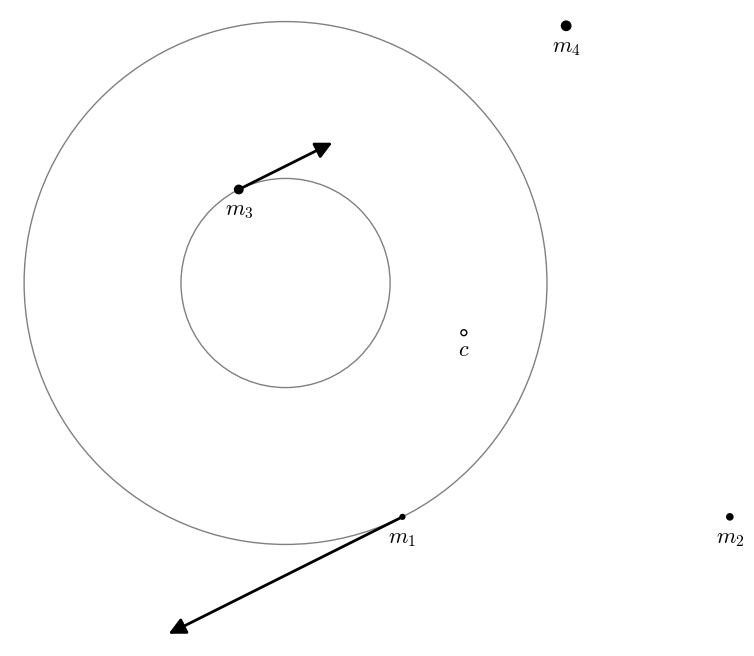}
\caption{A twist vector $v_{1,3}$ in a configuration of four masses, with $m_1=2$, $m_2=3$, $m_3=5$, and $m_4=6$. The system's center of mass $c$ is shown as an open circle.}

\label{Jacobitwist}
\end{figure}
\end{center}

An example of one of these vectors is shown in Figure (\ref{Jacobitwist}).

\begin{thm} The dimension of the subspace spanned by the (planar) $v_{i,j}$ is $2n-3$, unless all $n$ points are collinear, in which case the dimension is $n-1$.
\end{thm}
\begin{proof}
If all $N$ points are collinear, we can pick coordinates so that $x_i = 0$ for each $1\le i \le n$.  Then the $y_i$ components of the $v_{i,j}$ are all zero, and the $v_{i,j}$ are perpendicular to $\nabla c_x$, so the dimension of the span of the $v_{i,j}$ is at most $n-1$.  It is easy to see that $v_{1,j}$ are independent for $j>1$, so this subset is a basis for the span.

If the points are not all collinear, assume without loss of generality that $q_1$, $q_2$, and $q_3$ are not collinear.  By computing the determinant of the first few 3 by 3 minors of the matrix with rows $v_{1,2}$, $v_{1,3}$, $v_{2,3}$ it can be seen that these three vectors are independent.  Now we proceed by induction, assuming that we have a set $S_k$ of $2k-3$ independent vectors $v_{i,j}$ with $i < k+1$ and $j < k+1$.  We choose the next two vectors to be $v_{1,k+1}$ and one of $\{ v_{2,k+1}, v_{3,k+1}\}$.  If the points $q_1$, $q_2$, and $q_{k+1}$ are not collinear then we can choose $v_{2,k+1}$, otherwise we choose $v_{3,k+1}$.  These two new vectors are independent even in the projection onto the $x_{k+1}$ and $y_{k+1}$ coordinates, and thus are independent of the previous basis as well.   This process yields an independent subset of the $ \{v_{i,j}\} $ whose span is dimension $2 n-3$.   Since the $v_{i,j}$ are perpendicular to the three independent vectors $\nabla c_x$, $\nabla c_y$, and $\nabla I$, this is the maximum possible dimension of the span of any subset of the $\{v_{i,j}\}$.
\end{proof}

Applying the vector field $v_{i,j}$ to the function $f$ gives the well-known Laura-Andoyer equations for central configurations \cite{laura1905sulle, andoyer1906main}:

\begin{equation} \label{LAeqs}
\frac{ \partial f}{ \partial v_{i,j}} = \frac{ \partial (U/(A-2))}{ \partial v_{i,j}} = - m_i m_j \sum_{k} m_k ( R_{i,k} - R_{j,k}) \Lambda_{i,j,k} = 0 ,
\end{equation}
with the minor difference that our equations have an overall factor of $m_i m_j$, which is irrelevant if we assume the masses are nonzero.

\subsection{Use of twist vectors for numerical computations}

The vectors $v_{i,j}$ can be helpful in constructing numerical methods for finding central configurations and computing the Hessian of $f$.  As one example, we can use them in a gradient descent restricted to a constant value of the moment of inertia. Choosing a linearly independent set of the $v_{i,j}$, $V = \{V_1, V_2, \ldots V_m\}$, the flow generated by $\nabla_V(f)$ will usually converge to a local minimum of $f$ (after rescaling to a configuration with $MI = U$).  We have found that in some cases it can converge to a collinear configuration because of the lack of independence of the $V_i$ at such configurations.   

\begin{center}
\begin{figure}[h!t]
\includegraphics[width=3.25in]{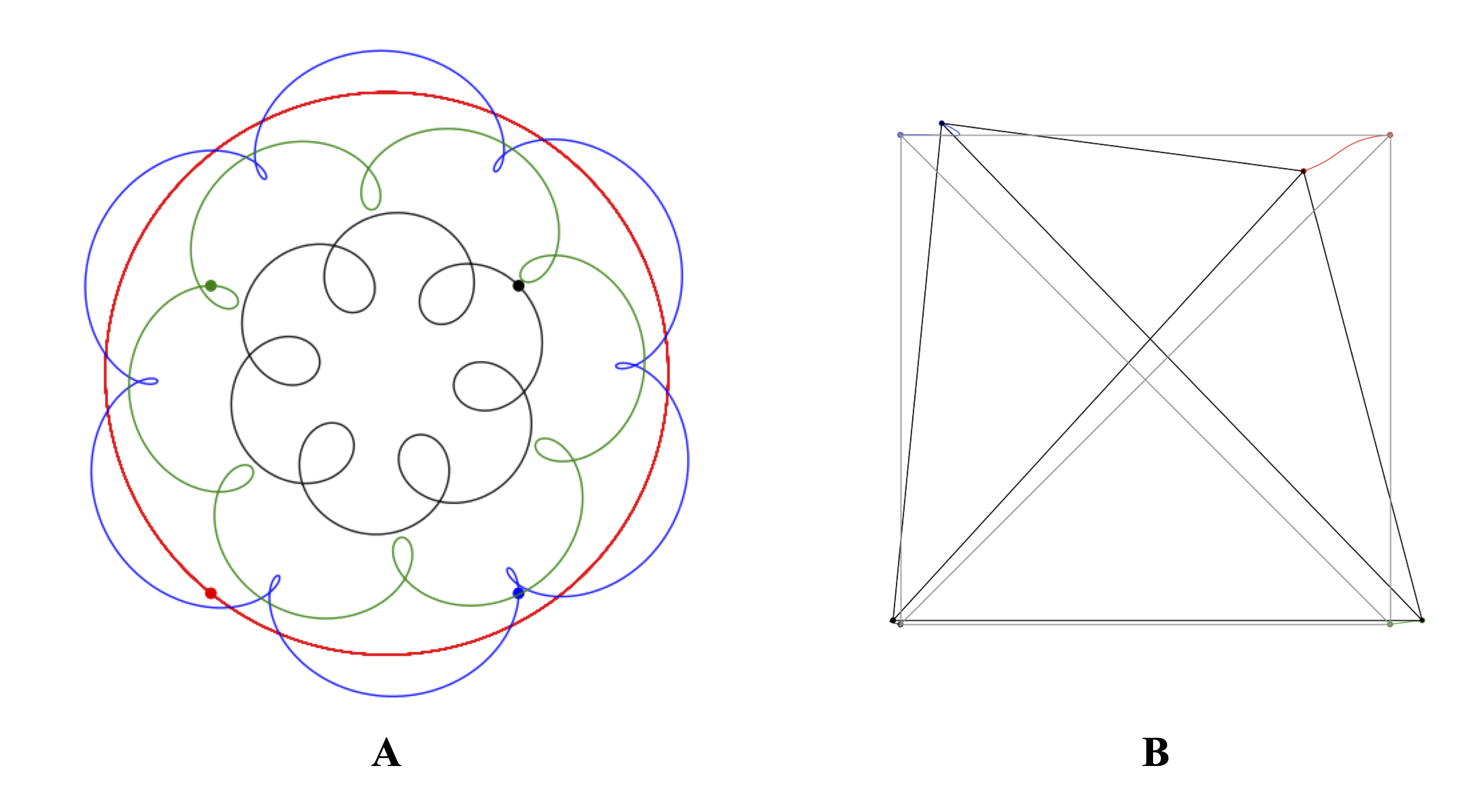}
\caption{Flows of a square configuration with masses 1,1,3, and 5 under two different twist vector fields: (A) a fixed linear combination of three elementary twists, and (B) a gradient descent to a convex central configuration}
\label{orb1123}
\end{figure}
\end{center}

Figure \ref{orb1123} shows a picture of the orbit of a square configuration with masses $m_1 = 1$,$m_2=1$,$m_3 = 2$, and $m_4 = 3$ using the fixed linear combination $v_{1,2} + v_{1,3} + v_{1,4}$.  

\subsection{The Hessian}

In Cartesian coordinates the components of the Hessian of $f$ are

$$ H_{x_i, x_k} = m_i m_k \left ( S_{i,k} - A x_{i,k}^2 T_{i,k} \right ), \ \ \ \   H_{y_i, y_k} = m_i m_k \left ( S_{i,k} - A y_{i,k}^2 T_{i,k} \right ) , $$

$$ H_{x_i, x_i} = - \sum_{k \neq i} H_{x_i, x_k},  \ \ \ \ H_{y_i, y_i} = - \sum_{k \neq i} H_{y_i, y_k} , $$

$$ H_{x_i, y_k} = - m_i m_j A T_{i,k} y_{i,k} x_{i,k} , $$

$$ H_{x_i, y_i} = - \sum_{k \neq i} H_{x_i, y_k} . $$

Using these formulae, we can compute the diagonal components of the Hessian with respect to the $v_{i,j}$:

\begin{equation} v_{i,j}^T H v_{i,j} = - m_i m_j \left ( (m_i + m_j)^2 r_{i,j}^2 S_{i,j} + \sum_{k\neq i,j} m_k \left [ r_{i,j}^2 ( m_j S_{i,k} + m_i S_{j,k} ) - A \Lambda_{i,j,k}^2 (m_j T_{i,k} + m_i T_{j,k}) \right ] \right ) .
\end{equation}

The second type of entry in the Hessian in these coordinates is of the form $v_{i,j}^T H v_{i,k}$ with $i$,$j$, and $k$ distinct:

\begin{equation}
\begin{split} 
 v_{i,j}^T H v_{i,k} = & m_i m_j m_k  \biggl (  q_{i,j} \cdot q_{i,k} \left [ m_i ( S_{k,j} - S_{i,j} - S_{i,k}) - \sum_{l \neq i}m_l S_{i,l} \right ]+ \\
 & - m_i A T_{j,k} \Lambda_{i,j,k}^2 - A \sum_{l \neq i,j,k} m_l T_{i,l} \Lambda_{i,j,l} \Lambda_{i,k, l} \biggr ) .
\end{split}
\end{equation}

The last type of entry is $v_{i,j}^T H v_{k,l}$ where all of $i,j,k,l$ are distinct from each other:

\begin{equation}
\begin{split}  v_{i,j}^T H v_{k,l} = & m_i m_j m_k m_l \biggl ( -A ( T_{i,k} \Lambda_{i,k,l} \Lambda_{i,k,j} + T_{j,k} \Lambda_{j,k,l} \Lambda_{j,k,i} +  \\ & T_{i,l} \Lambda_{i,l,k} \Lambda_{i,l,j} +  T_{j,l} \Lambda_{j,l,i} \Lambda_{j,l,k} ) + q_{i,j} \cdot q_{k,l} (R_{i,k} + R_{j,l} - R_{j,k} - R_{i,l})\biggr ) .
\end{split}
\end{equation}

Alternatively, if we just want to compute the signature of the Hessian at critical points of $f$, we can rescale the $v_{i,j}$ to 

$$ \tilde{v}_{i,j} = \frac{v_{i,j}}{m_i m_j r_{i,j}} . $$
Then using the notation $\mu_{i,j} = \frac{m_i}{m_j}$, $\theta_{i,j,k}$ for the angle between vectors $q_{i,j}$ and $q_{i,k}$, and $\theta_{i,j;k,l}$ for the angle between vectors $q_{i,j}$ and $q_{k,l}$, the Hessian components become

\begin{equation} \tilde{v}_{i,j}^T H \tilde{v}_{i,j} =  -(2 + \mu_{i,j} + \mu_{j,i}) S_{i,j} + \sum_{k\neq i,j} \left [  -\mu_{k,i} S_{i,k} - \mu_{k,j} S_{j,k}  + A  (\mu_{k,i} \sin^2 \theta_{j,i,k} R_{i,k} + \mu_{k,j}  \sin^2 \theta_{i,j,k}  R_{j,k})  \right ] ,
\end{equation}

\begin{equation}
\begin{split} 
\tilde{v}_{i,j}^T H \tilde{v}_{i,k} = & \cos \theta_{j,i,k} \left [ S_{k,j} - S_{i,j} - S_{i,k}  - \sum_{l \neq i} \mu_{l,i} S_{i,l} \right ] \\
 & - A R_{j,k} \sin \theta_{i,j,k} \sin \theta_{i,k,j} - A \sum_{l \neq i,j,k} \mu_{l,i} R_{i,l} \sin \theta_{j,i,l} \sin \theta_{k,i, l} ,
\end{split}
\end{equation}

\begin{equation}
\begin{split} \tilde{v}_{i,j}^T H \tilde{v}_{k,l} = &   -A ( R_{i,k} \sin \theta_{i,k,l} \sin \theta_{j,i,k} + R_{j,k} \sin \theta_{j,k,l} \sin \theta_{i,j,k}  +  R_{i,l} \sin \theta_{i,l,k} \sin \theta_{j,i,l}  \\ & +  R_{j,l} \sin \theta_{j,l,k} \sin \theta_{i,j,l} ) + \cos \theta_{i,j;k,l} (R_{i,k} + R_{j,l} - R_{j,k} - R_{i,l}).
\end{split}
\end{equation}

\section{Application to the 3-body problem}

We will give a brief overview of the application of our equations to the 3-body case, which has been so well studied that no new results will be obtained, but it is perhaps useful to see them in this context.

For three point masses the sum in each of Laura-Andoyer equations (\ref{LAeqs}) reduces to a single term:

$$ \frac{ \partial f}{ \partial v_{i,j}} = - m_i m_j m_k ( R_{i,k} - R_{j,k}) \Lambda_{i,j,k} = 0 $$
in which the $i$, $j$, and $k$ are distinct.  If $\Lambda_{i,j,k}$ is nonzero then $R_{i,k} = R_{j,k}$ and the points must be equidistant, so the only non-collinear central configurations are the Lagrange equilateral triangles. Any three collinear points satisfy these equations, but the $v_{i,j}$ are not independent and some additional equations must be used to recover the Euler solutions.  

For the equilateral triangle with our choice of $f$, the distances $r_{i,j} = 1$, and so $S_{i,j} = 0$, $R_{i,j} = T_{i,j} = 1$, all of the angles are $\pi/3$ so $\sin^2 \theta_{i,j,k} = 3/4$.  We can compute the Hessian components:

$$ \tilde{v}_{i,j}^T H \tilde{v}_{i,j} = \frac{3 A}{4} (\mu_{k,i} + \mu_{k,j}) , $$

$$ \tilde{v}_{i,j}^T H \tilde{v}_{i,k} =  - \frac{3 A}{4} , $$
so in coordinates $(\tilde{v}_{1,2},\tilde{v}_{1,3},\tilde{v}_{2,3})$ the Hessian matrix is

$$ \frac{3 A}{4} \left ( \begin{array}{ccc}  \mu_{3,1} + \mu_{3,2} & -1 & -1 \\ -1 & \mu_{2,1} + \mu_{2,3} & -1 \\  -1 & -1 &\mu_{1,2} + \mu_{1,3}   \end{array} \right ).$$

This matrix has one zero eigenvalue corresponding to the rotational symmetry of the space, and two positive eigenvalues.  To prove positivity of the nonzero eigenvalues we first note that the matrix is symmetric, so its eigenvalues are real.  Its characteristic polynomial factors as $\lambda (a_0 + a_1 \lambda  + a_2 \lambda)$, in which the coefficients $a_i$ are functions of the mass ratios $\mu_{j,k}$.  The reality of the eigenvalues means that the discriminant $a_1^2 - 4 a_0 a_2 \ge 0$, and it is straightforward to compute that the functions $a_0$ and $a_2$ have the same sign, so $ a_1 > \sqrt{a_1^2 - 4 a_0 a_2}$.  It is also straightforward to compute that 

$$ a_1 = -( \mu_{1,2} + \mu_{1,3} + \mu_{2,3} + \mu_{2,1} + \mu_{3,1} + \mu_{3,2})  < 0 $$
so the eigenvalues $-a_1/2 \pm \sqrt{a_1^2 - 4 a_0 a_2}/2$ are positive, and the Lagrange configurations are minima of $f$.

\section{Application to the 4-body problem}

We will illustrate the utility of our new equations for the Hessian by computing the Morse index of the four-body central configurations which are concave isosceles triangles (sometimes called concave kites).  This includes the case of equilateral triangles with one arbitrary nonnegative central mass and three equal outer masses.  This is an interesting special case since it has long been known to possess a degenerate Hessian for particular ratios of the central to outer masses \cite{palmore73,meyer1988bifurcationsB,meyer1988bifurcations,fbdbif2017}.  

There are several bifurcations in the index of the Hessian, which correspond to some of the changes in number of concave central configurations seen for example in \cite{simo_relative_1978}.  

Our coordinates for the four points will be $q_1 = (-x_2,0)$, $q_2 = (x_2, 0)$, $q_3 = (0,y_3)$, and $q_4 = (0,y_4)$.  We will assume that $x_2$ is positive and $y_3 > y_4$.  With elementary geometry we can compute: 

$$ \Lambda_{1,2,3} = 2 x_2 y_3, \ \ \ \Lambda_{1,2,4} = 2 x_2 y_4, \ \ \  \Lambda_{1,3,4} = - \Lambda_{2,3,4} = x_2 (y_3 -  y_4), $$

$$\sin \theta_{1,2,3}  = \frac{y_3}{r_{1,3}}, \ \ \ \sin \theta_{1,2,4}  = \frac{y_4}{r_{1,4}} , \ \ \ \sin \theta_{1,3,2}  = \frac{r_{1,2} y_3}{r_{1,3}^2}, \ \ \ \sin \theta_{1,4,2} = \frac{r_{1,2}y_4}{r_{1,4}^2}, $$

$$ \sin \theta_{1,3,4}  = \frac{r_{1,2}}{2 r_{1,3}}, \ \ \ \ \sin \theta_{3,1,4}  = \frac{r_{3,4} r_{1,2}}{2 r_{1,3} r_{1,4}},\ \ \ \ \sin \theta_{1,4,3} = \frac{r_{1,2}}{2 r_{1,4}} ,$$

$$ \cos \theta_{1,3,2} = 1 - \frac{r_{1,2}^2}{2 r_{1,3}^2}, \ \ \ \cos \theta_{1,3;2,4} = \frac{y_3 y_4 - x_2^2}{r_{1,3} r_{1,4}}. $$

One of the critical point equations (\ref{LAeqs}) implies that $m_1 = m_2$, and the others which are nontrivial simplify to

$$ \frac{1}{m_1 m_3 }\frac{\partial f}{\partial v_{1,3}} = - 2 m_1 (R_{1,2} - R_{1,3}) x_2 y_3 + m_4 (R_{1,4} - R_{3,4}) x_2 r_{3,4} = 0 , $$

$$ \frac{1}{m_1 m_4 }\frac{\partial f}{\partial v_{1,4}} = - 2 m_1 (R_{1,2} - R_{1,4}) x_2 y_4 - m_3 (R_{1,3} - R_{3,4}) x_2 r_{3,4} = 0 , $$
so the mass ratios are

\begin{equation} \label{mukite1} \mu_{3,1} = \frac{m_3}{m_1} = \frac{2 y_4}{r_{3,4}} \frac{(R_{1,4} - R_{1,2})}{(R_{3,4} - R_{1,3})} , \end{equation}

\begin{equation}  \label{mukite2} \mu_{4,1} = \frac{m_4}{m_1} = \frac{2 y_3}{r_{3,4}} \frac{(R_{1,2} - R_{1,3})}{(R_{1,4} - R_{3,4})} . \end{equation}

Besides the equilateral case, there are two classes of concave kites with non-negative masses.  If $\frac{y_3}{x_2} > \sqrt{3}$, then $r_{1,3} = r_{2,3} > r_{1,2}$ and $\frac{y_3}{2} - \frac{x_2^2}{2 y_3} < y_4 < \frac{x_2 \sqrt{3}}{2}$.  Alternatively if $y_3 < \frac{\sqrt{3}}{2}$ then for non-negative masses we must have $0 < y_4 < \frac{y_3}{2} - \frac{x_2^2}{2 y_3}$.  When $y_4 = \frac{y_3}{2} - \frac{x_2^2}{2 y_3}$, the fourth point is at the circumcenter of the outer triangle, and $r_{1,4} = r_{2,4} = r_{3,4}$.

\begin{figure}[h!t]
\includegraphics[width=4in]{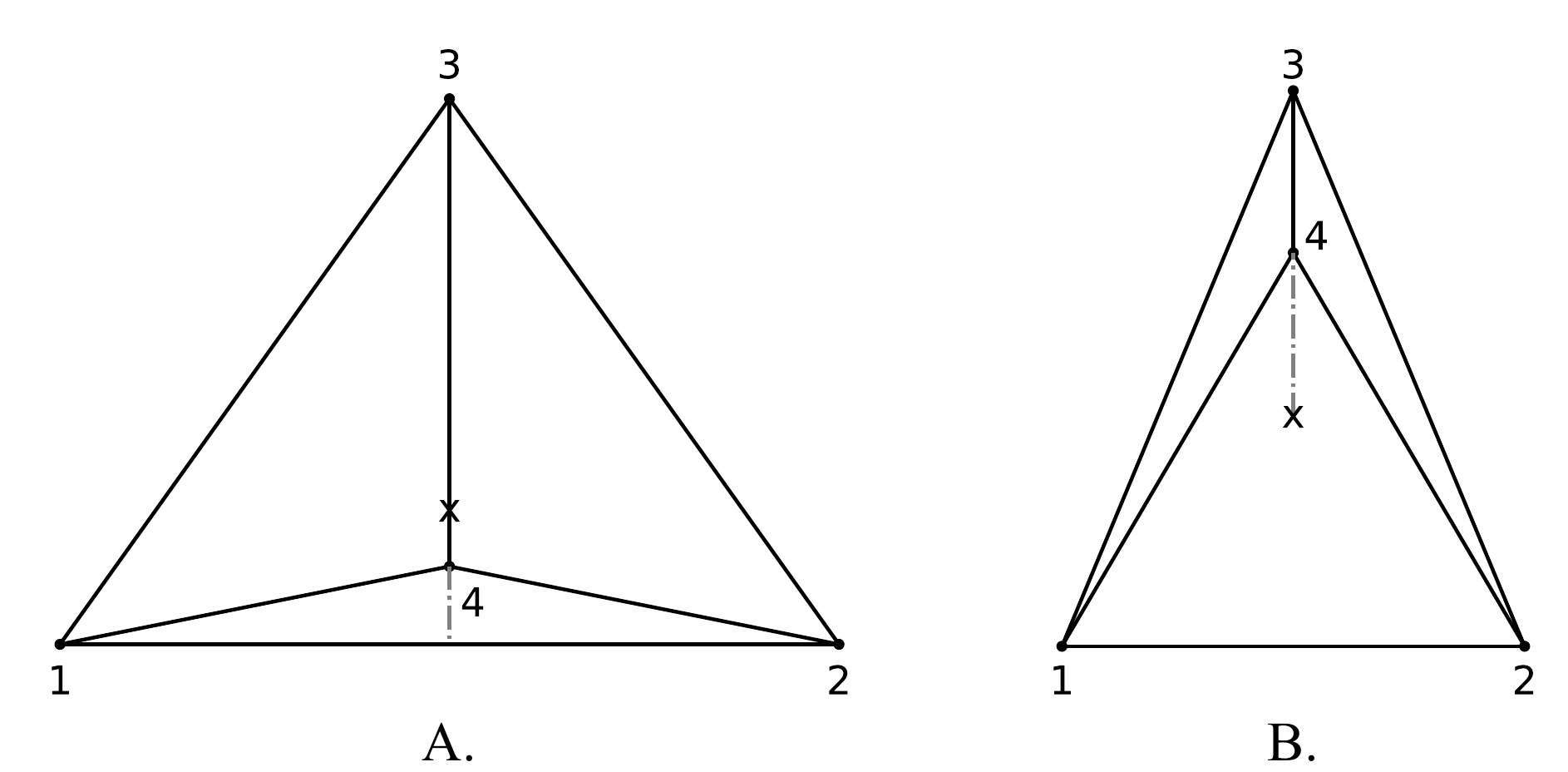}
\caption{Two examples of concave kites.  A. An example `wide' kite with $r_{1,3}=r_{2,3} < r_{1,2}$; the fourth point must be below the circumcenter (marked with an `x') for positive mass configurations. B. An example `tall' kite with $r_{1,3}=r_{2,3} > r_{1,2}$; the fourth point must be above the circumcenter.}
\label{ExampleKites}
\end{figure}

The Dziobek equations \cite{dziobek_ueber_1900, moeckel_generic_2001}:

$$ S_{1,2} S_{3,4} = S_{1,3} S_{2,4} = S_{1,4} S_{2,3} $$
in the case of the isosceles kites specialize to just one equation, 

\begin{equation} \label{DzIso}
S_{1,2} S_{3,4} = S_{1,3} S_{1,4}.
\end{equation}

This Dziobek equation can be used to find $x_2 = r_{1,2}/2$ given the shape variables $z_3 = \frac{y_3}{x_2}$ and $z_4 = \frac{y_4}{x_2}$:

\begin{equation} \label{Dzsub}
x_2 = \frac{1}{2} \left (\frac{(z_3 - z_4)^{-A} - (1+z_3^2)^{-A}(1+z_4^2)^{-A}}{1 + (z_3 - z_4)^{-A} -(1+z_3^2)^{-A}-(1+z_4^2)^{-A}} \right )^{1/A} .
\end{equation}

We can block-diagonalize the Hessian for these kites by using the basis

$$ \beta = [ \tilde{v}_{a,3} = \frac{1}{\sqrt{2}} (\tilde{v}_{1,3} + v_{2,3}), \ \ \  \tilde{v}_{a,4} = \frac{1}{\sqrt{2}} (v_{1,4} + v_{2,4}), \ \ \ \  \tilde{v}_{3,4},$$
$$\tilde{v}_{s,3} = \frac{1}{\sqrt{2}} (\tilde{v}_{1,3} - \tilde{v}_{2,3}), \ \ \ \tilde{v}_{s,4} = \frac{1}{\sqrt{2}} (v_{1,4} - v_{2,4}) ] ,$$
which results in blocks of dimensions $3 \times 3$ and $2 \times 2$, which we will denote by $H_a$ and $H_s$ respectively.

The $H_a$ block always contains the eigenvector of eigenvalue 0 which is tangent to the uniform rotation about the center of mass.  While it is possible to symbolically calculate this eigenvector in these coordinates, it is quite complicated, and it does not seem productive to calculate the resulting residual $2 \times 2$ submatrix explicitly.

The $2 \times 2$ block $H_s$ is easier to analyze, and it contains an interesting bifurcation.   We will focus on the configurations that are close to the coorbital case, in which the fourth point is close to the circumcenter of the outer three points.  The first component of $H_s$ is

$$ (H_s)_{1,1} = \frac{1}{2} \left ( \tilde{v}_{1,3}^T H \tilde{v}_{1,3} + \tilde{v}_{2,3}^T H \tilde{v}_{2,3} - 2 \tilde{v}_{1,3}^T H \tilde{v}_{2,3} \right ) $$
$$ = \tilde{v}_{1,3}^T H \tilde{v}_{1,3}  - \tilde{v}_{1,3}^T H \tilde{v}_{2,3} $$

$$ =  -(2 + \mu_{1,3} + \mu_{3,1}) S_{1,3} + \sum_{k\neq 1,3} \left [  -\mu_{k,1} S_{1,k} - \mu_{k,3} S_{3,k}  + A  (\mu_{k,1} \sin^2 \theta_{3,1,k} R_{1,k} + \mu_{k,3}  \sin^2 \theta_{1,3,k}  R_{3,k})  \right ]   $$

$$ - \left(   \cos \theta_{1,3,2} \left [ S_{1,2} - S_{1,3} - S_{2,3}  - \sum_{l \neq 3} \mu_{l,3} S_{3,l} \right ] - A R_{1,2} \sin \theta_{3,1,2} \sin \theta_{3,2,1} - A \mu_{4,3} R_{3,4} \sin \theta_{1,3,4} \sin \theta_{2,3,4}  \right ) $$

$$ = -(2 + 2 \mu_{1,3} + \mu_{3,1}) S_{1,3} - S_{1,2} - \mu_{4,1}S_{1,4} - \mu_{4,3} S_{3,4} - \left ( 1 - \frac{r_{1,2}^2}{2 r_{1,3}^2} \right ) \left ( S_{1,2} - 2 S_{1,3} - 2 \mu_{1,3} S_{1,3} - \mu_{4,3} S_{3,4} \right )   $$

$$ + A \left ( \frac{2 y_3^2}{r_{1,3}^2} R_{1,2} + \mu_{2,3}\frac{r_{1,2}^2 y_3^2}{r_{1,3}^2} R_{1,3} + \mu_{4,1} \frac{r_{1,2}^2 r_{3,4}^2 }{4 r_{1,3}^2 r_{1,4}^2 } R_{1,4} + 2 \mu_{4,3} \frac{r_{1,2}^2}{4 r_{1,3}^2} R_{3,4}\right ) . $$

Similarly

$$(H_s)_{2,2} = \tilde{v}_{1,4}^T H \tilde{v}_{1,4} - \tilde{v}_{1,4}^T H \tilde{v}_{2,4} $$

$$ =  -(2 + \mu_{1,4} + \mu_{4,1}) S_{1,4} + \sum_{k\neq 1,4} \left [  -\mu_{k,1} S_{1,k} - \mu_{k,4} S_{4,k}  + A  (\mu_{k,1} \sin^2 \theta_{4,1,k} R_{1,k} + \mu_{k,4}  \sin^2 \theta_{1,4,k}  R_{4,k})  \right ]   $$

$$ - \left(   \cos \theta_{1,4,2} \left [ S_{1,2} - S_{4,1} - S_{4,2}  - \sum_{l \neq 4} \mu_{l,4} S_{4,l} \right ] - A R_{1,2} \sin \theta_{4,1,2} \sin \theta_{4,2,1} - A \mu_{3,4} R_{3,4} \sin \theta_{1,4,3} \sin \theta_{2,4,3}  \right ) $$

$$ = -(2 + 2 \mu_{1,4} + \mu_{4,1}) S_{1,4} - S_{1,2} - \mu_{3,1}S_{1,3} - \mu_{3,4} S_{3,4}- \left ( 1 - \frac{r_{1,2}^2}{2 r_{1,4}^2} \right ) \left ( S_{1,2} - 2 S_{1,4} - 2 \mu_{1,4} S_{1,4} - \mu_{3,4} S_{3,4} \right )   $$

$$ + A \left ( \frac{2 y_4^2}{r_{1,4}^2} R_{1,2} + \mu_{2,4}\frac{r_{1,2}^2 y_4^2}{r_{1,4}^2} R_{1,4} + \mu_{3,1} \frac{r_{1,2}^2 r_{3,4}^2 }{4 r_{1,3}^2 r_{1,4}^2 } R_{1,3} + 2 \mu_{3,4} \frac{r_{1,2}^2}{4 r_{1,4}^2} R_{3,4}\right ) , $$

and 
$$(H_s)_{1,2} = (H_s)_{2,1} = \tilde{v}_{1,3}^T H \tilde{v}_{1,4} + \tilde{v}_{1,3}^T H \tilde{v}_{2,4} , $$

in which 
\begin{equation*}
  \begin{aligned}
\tilde{v}_{1,3}^T H \tilde{v}_{2,4} & = \cos{\theta_{1,3;2,4}} ( R_{1,2} + R_{3,4} - R_{1,3} - R_{1,4}) - A (R_{1,2} \sin{\theta_{1,2,4}} \sin{\theta_{3,1,2}}  \\
+ & R_{1,3} \sin{\theta_{3,1,4}} \sin{\theta_{1,3,2}} + R_{1,4} \sin{\theta_{1,4,2}} \sin{\theta_{3,1,4}} + R_{3,4} \sin{\theta_{3,4,2}} \sin{\theta_{1,3,4} }) \\
& =  \frac{1}{r_{1,3}r_{1,4}} \left [ (y_3 y_4 - x_2^2) (R_{1,2} +R_{3,4} - R_{1,3} - R_{1,4} ) \right . \\
  & \left .  -  A \left (R_{1,2} y_3 y_4 + R_{1,3} \frac{r_{1,2}^2 r_{3,4} y_3}{2 r_{1,3}^2} + R_{1,4} \frac{r_{1,2}^2 r_{3,4} y_4}{2 r_{1,4}^2} + R_{3,4} \frac{r_{1,2}^2}{4} \right )\right ] ,
   \end{aligned}
 \end{equation*}
 
 and
 % i = 1, j = 3, k = 4
 \begin{equation*}
  \begin{aligned}
  \tilde{v}_{1,3}^T H \tilde{v}_{1,4} & = \cos \theta_{3,1,4} \left ( S_{3,4} - S_{1,3} - S_{1,4} - \mu_{2,1} S_{1,2} - \mu_{3,1} S_{1,3} - \mu_{4,1} S_{1,4} \right ) \\
  & - A R_{3,4} \sin \theta_{1,3,4} \sin \theta_{1,4,3} - A \mu_{2,1} R_{1,2} \sin \theta_{3,1,2} \sin \theta_{1,2,4} \\
  & =  \frac{1}{4 r_{1,3} r_{1,4}} \left \{ (r_{1,2}^2+ 4 y_3 y_4) \Bigl[ S_{3,4} - (1 + \mu_{3,1}) S_{1,3} - (1 + \mu_{4,1})S_{1,4} - S_{1,2} \Bigr ] \right . \\
  & \left . - A ( R_{3,4} r_{1,2}^2  - 4 R_{1,2} y_3 y_4 ) \right \} .
    \end{aligned}
 \end{equation*}
 
The eigenvalues of these Hessians contain several bifurcations depending the shape of the isosceles configuration and the exponent parameter $A$, which are difficult to analytically characterize in general.  For the Newtonian case we investigated the Hessian eigenvalues more thoroughly numerically.  Most of these numerical conclusions can be refined to rigorous results using interval arithmetic, since there are only two shape parameters, and our block diagonalization of the Hessian greatly reduces the computations. 

\begin{figure}[h!b] 
\includegraphics[height=4in]{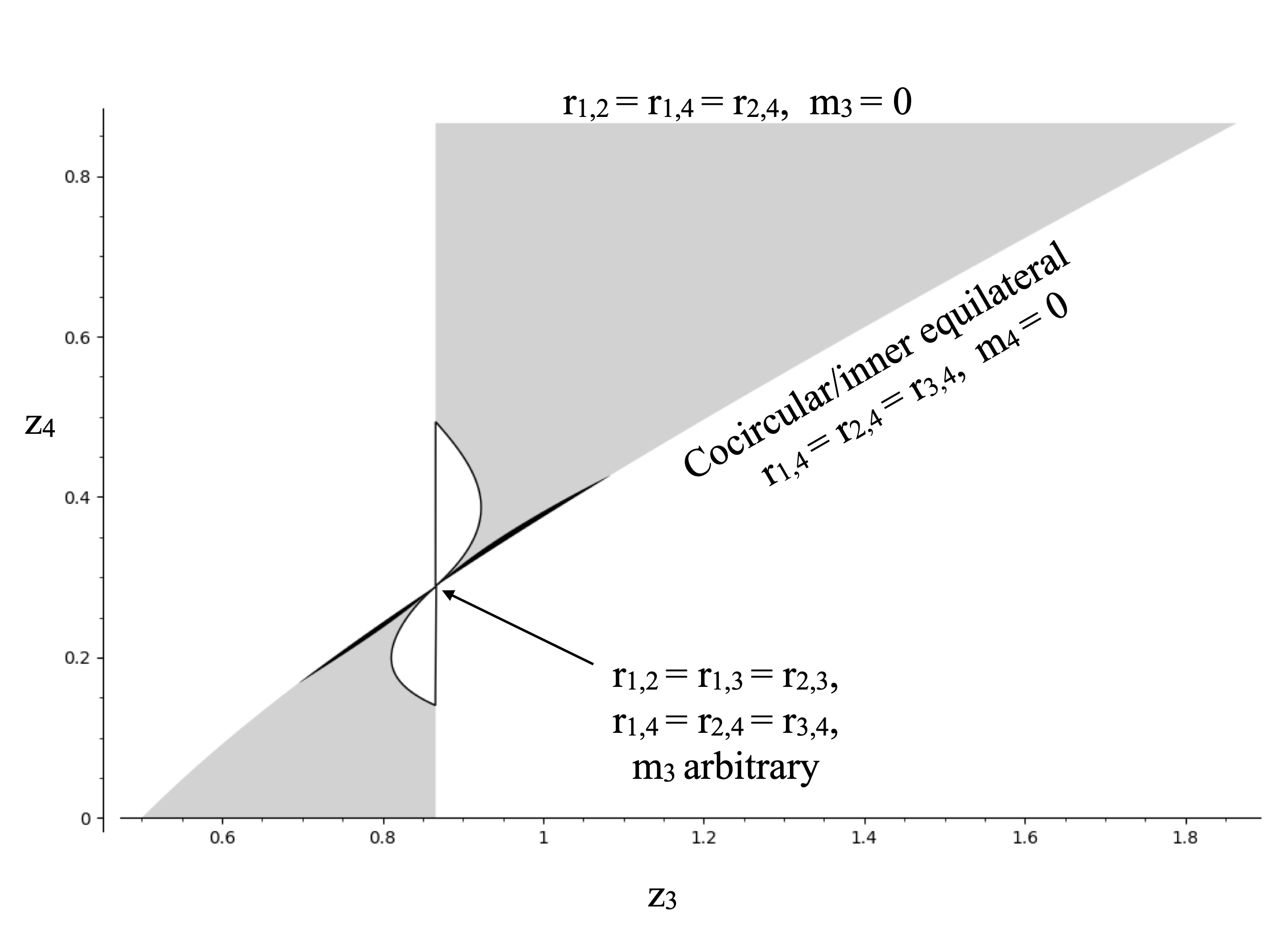}
\caption{Isosceles kite central configurations for the Newtonian potental with nonnegative masses parameterized by $z_3$ and $z_4$.  Configurations with a positive-definite Hessian are shaded in white; Hessians of index 3 are shaded in grey, and Hessians of index 2 are black.}
\label{FigIsoKites}
\end{figure}

In Figure \ref{FigIsoKites} we illustrate the concave isosceles configurations parameterized by $z_3$ and $z_4$.  

When  $z_4 = y_c$ the fourth point is at the circumcenter of the outer triangle formed by the first three points, so $r_{1,3} = r_{2,3} = r_{3,4}$.  To have a central configuration with positive masses we must have $z_4 > y_c$.  As $z_4$ approaches $y_c$ we obtain in the limit the $1+3$ coorbital central configurations, which have been extensively characterized already  \cite{RS2004, corbera2015bifurcation, hampton_deng_coorbital}.

From equations (\ref{mukite1}) and (\ref{mukite2}) we see that for a fixed outer triangle shape as $z_4$ approaches $y_c$ the ratio of $m_4$ to the outer masses becomes infinite.  Some terms in the Hessian that contain $\mu_{4,1}$ or $\mu_{4,3}$ will then dominate the others.  There is some subtlety in products such as $\mu_{4,1} S_{1,4}$ however, since the inner distances become close to $1$ in this limit,  $S_{1,4}$ and $S_{3,4}$ approach 0.  We can handle these quantities by using the critical point relation $U = MI$ described earlier, and the Dziobek equation (\ref{DzIso}).  From this relation we can see that $R_{1,4}$ and $R_{3,4}$ must approach $1$ as $z_4$ approaches $y_c$ if the outer shape is fixed.  

If we divide the relation $U - MI = 0$ by $m_1^2$, for the isosceles central configurations we have

$$ r_{1,2}^2 S_{1,2} + 2 \mu_{3,1} r_{1,3}^2 S_{1,3} + 2 \mu_{4,1} r_{1,4}^2 S_{1,4}  + \mu_{4,1} \mu_{3,1} r_{3,4}^2 S_{3,4} = 0. $$

Eliminating $S_{3,4}$ with the Dziobek equation we can express $\mu_{4,1} S_{1,4}$ in terms of quantities with more obvious limits:

$$ \mu_{4,1} S_{1,4} = - \frac{r_{1,2}^2 S_{1,2} + 2 \mu_{3,1} r_{1,3}^2 S_{1,3}} {2r_{1,4}^2 + \mu_{3,1}r_{3,4}^2 S_{1,3}/S_{1,2}}. $$

The $H_s$ submatrix of the Hessian is responsible for the upper ($z_3 > \sqrt{3}/2$) black to grey bifurcation in Figure \ref{FigIsoKites}, where the Hessian changes from index 2 to index 1, and the lower white to grey bifurcation ($z_3 < \sqrt{3}/2$) in which the Hessian changes from index 0 to index 1.  In the index 2 configurations, one negative eigenvalue comes from $H_a$ and one from $H_s$, with the boundary determined by the change of a positive to a negative eigenvalue in $H_s$ or $H_a$ for the upper or lower change, respectively.  Similarly, the white shaded region corresponds to configurations with index 0, the boundary of which is defined by a zero eigenvalue of $H_a$ for the upper region and $H_s$ for the lower region.  These switches between the responsible submatrix were surprising to us.
=
We can verify the picture in Figure \ref{FigIsoKites} with interval arithmetic to a resolution which is at least an order of magnitude better than in the figure.  It is difficult to completely prove the indicated structure of the bifurcations, especially near the boundaries where some of the masses approach zero, and the most symmetric point where  $r_{1,2} = r_{1,3}=r_{2,3}$ and $r_{1,4} = r_{2,4}=r_{3,4}$.   Some changes in variables similar to equation (\ref{Dzsub}) could probably be used near the boundaries to fix these issues, but there the problem reduces to the restricted problem, which has been previously well characterized \cite{Lindlow22, pedersen1944librationspunkte, pedersen1952stabilitatsuntersuchungen, Gannaway81, Arens, leandro_central_2006, barros2014bifurcations}.   Near the symmetric point, where the mass $m_4$ is arbitary, a perturbative analysis has been done for some aspects of the central configurations \cite{meyer1988bifurcationsB}.

It would be interesting to see if these bifurcation curves could be computed exactly, for example with Gr\"obner bases, which seems quite possible; we did not attempt that.

\section{Conclusion and future directions}

To the best of our knowledge the framework of the twist vectors introduced here has not been described before, although because of the relative simplicity of the idea it is very possible that it occurs elsewhere in the classical mechanics literature.  

Although our equations are not very simple, in a standard Cartesian basis the Hessian of the function $f$ is much more complicated and harder to interpret in terms of geometric quantities.  We hope the approach presented here can be used to improve results on the finiteness, enumeration, and stability analysis of relative equilibria and central configurations in the $N$-body family of problems.

\bibliography{../../CelMechEtc}{}
\bibliographystyle{amsplain}

\end{document}